\documentclass{amsart}
\usepackage[dvips]{graphics}
\usepackage{comment}
\usepackage{pstricks}
\usepackage{pst-node}
\newtheorem{theorem}{Theorem}[section]
\newtheorem{propo}[theorem]{Proposition}
\newtheorem{lemma}[theorem]{Lemma}
\newtheorem{coro}[theorem]{Corollary}
\definecolor{brightgreen}{rgb}{0.85,1,0.85}
\definecolor{brightyellow}{rgb}{1,1,0.7}
\def\C{\mathbb{C}}
\newcommand{\fb}{ \psframebox[linearc=0.2,cornersize=absolute, framesep=6pt,fillstyle=solid,fillcolor=brightgreen]}
\newcommand{\gfb}{ \psframebox[linearc=0.2,cornersize=absolute, framesep=6pt,fillstyle=solid,fillcolor=brightyellow]}

\title{Selfsimilarity in the Birkhoff sum of the cotangent function}
\date{March 22, 2012}
\author{Oliver Knill}
\email{Knill@math.harvard.edu}
\address{
        Department of Mathematics \\
        Harvard University \\
        Cambridge, MA, 02138\\
        }

\begin{document}
\maketitle

\begin{abstract}
We prove that the Birkhoff sum $S_n(\alpha) = \frac{1}{n} \sum_{k=1}^{n} g(k \alpha)$
of $g = \cot(\pi x)$ and with golden ratio $\alpha$ converges in the sense that the sequence of
functions $s_n(x) = S_{[x q_{2n}]}/q_{2n}$ with Fibonacci $q_n$
converges to a self similar limiting function $s(x)$ on $[0,1]$ which can be computed analytically.
While for any continuous function $g$, the Birkhoff limiting function is $s(x)=M x$ by Birkhoff's ergodic
theorem, we get so examples of random variables $X_n$, where the limiting function of
$S_{[x n]}/n \to s(x)$ exists along subsequences for one initial point and is nontrivial.
Hardy and Littlewood have studied the Birkhoff sum for 
$g'(x)= -\pi \csc^2(\pi x)$ and shown that $\frac{1}{n^2} \sum_{k=1}^{[x n]} g'(k \alpha)$
stays bounded. Sinai and Ulcigrai have found a limiting distribution for $g(x)$
if both the rotation number $\alpha$ and the initial point $\theta$ are integrated over.  
We fix the golden ratio $\alpha$, start with fixed $\theta=0$ and show that
the rescaled random walk converges along subsequences.  
\end{abstract}

\section{Introduction} 

We study the sum $S_n = \sum_{k=1}^{n} g(k \alpha)$ with unbounded $g(x)$ and irrational $\alpha$. 
If the function $g$ is not in $L^1(T^1)$, there is no ergodic theorem to compute the limit of the 
Birkhoff averages $S_n(\alpha)/n$.
We will see however that there is a convergence in a rather particular way if $\alpha$ is 
the golden mean and $g$ has a single pole like for $g(x) = \cot(\pi x)$.
We have studied this particular Birkhoff sum for the antiderivative $G(x)=\log|1-e^{2\pi i x}| \frac{1}{\pi}$ of $g$ 
in \cite{KTp} with the help of Birkhoff functions $s_m(x) =\frac{S_{[m x]}}{m}$ on $[0,1]$. 
The function $s_m(x)$ encodes the random walk $\{S_k\}_{k=1}^m$. For identically distributed $L^1$ random variables, 
Birkhoff's ergodic theorem $X_k(\theta)=g(\theta+k \alpha)$ assures that $s_m(x) \to M x$, 
pointwise for almost all $\theta$, where 
$M = {\rm E}[X]=\int_0^1 g(\theta) \; d\theta$ is the expectation. We call the limit $s(x) = Mx$ the {\bf Birkhoff limiting
function} of $g$.  \\
Our main result is to that that for the random variables $X_k(\theta) = \cot(\pi (\theta + k \alpha))$ 
and particular starting point $\theta=0$, the functions $s_m(x)$ converge pointwise to a selfsimilar Birkhoff 
limiting function $s(x)$ along subsequences $q_n$, where
$p_n/q_n$ denote partial fractions of $\alpha$ so that $q_n=p_{n+1}$ are the Fibonacci numbers:

\begin{theorem}
\label{main}
The {\bf Birkhoff limiting function} 
$$  s(x) = \lim_{n \to \infty} \frac{S_{[q_n x]}}{q_n} \; , $$
for $g(x) = \cot(\pi x)$ exists pointwise along odd and even subsequences. 
The graph of $s(x) = \lim_{n \to \infty} s_{2n}(x)$ 
satisfies $s(\alpha x) = -\alpha s(x)$ and is continuous from the right. 
\end{theorem}

\begin{figure}
\scalebox{0.9}{\includegraphics{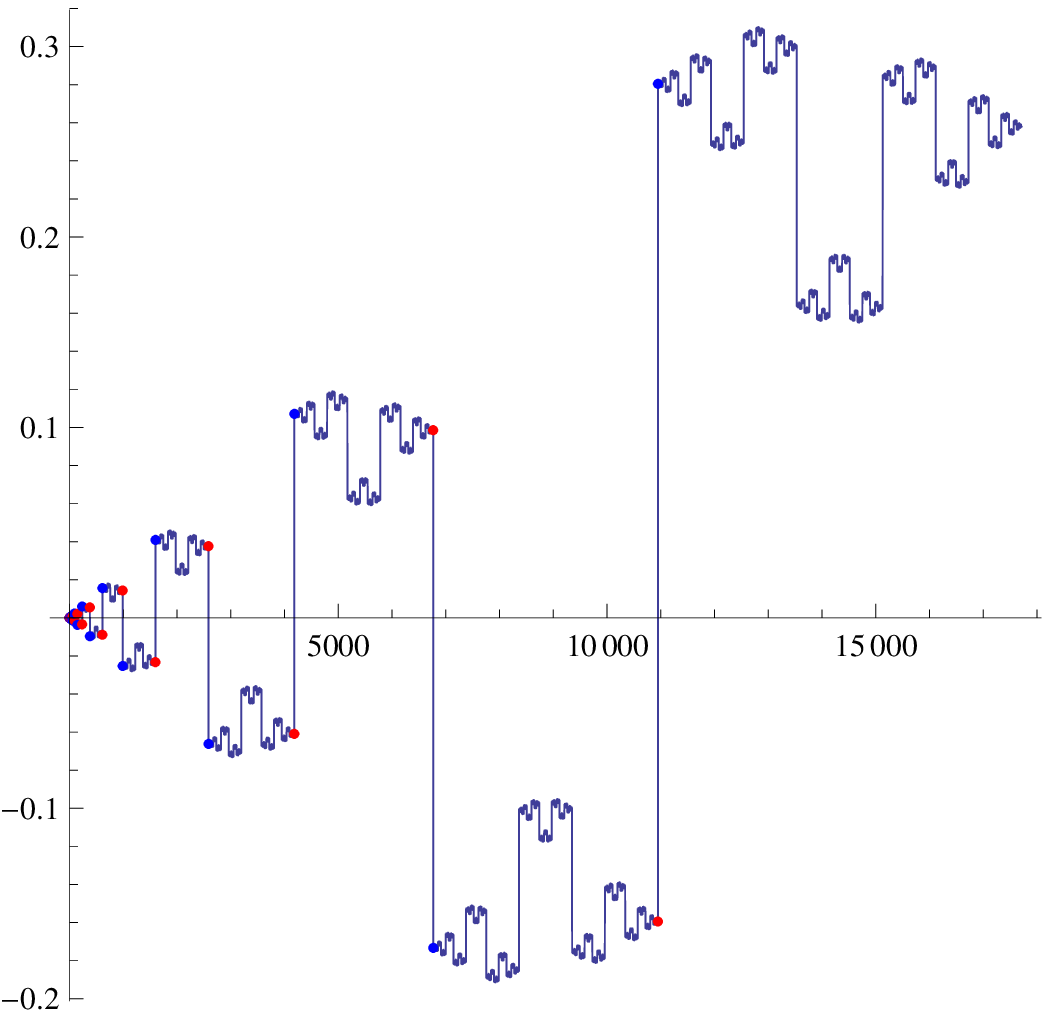}} \\
\caption{
The {\bf golden graph} is the graph of the Birkhoff limiting function
$$  s(x) =   \lim_{n \to \infty} S_{[q_{2n} x]}/q_{2n} 
         = - \lim_{n \to \infty} S_{[q_{2n+1} x]}/q_{2n+1}  \; . $$
This limit is more interesting than for continuous functions $g$, where by Birkhoff's ergodic theorem, the Birkhoff
limiting function is a linear function $s(x)=Mx$.
}
\label{goldengraph}
\end{figure}

We call the graph of $s$ the {\bf golden graph}. It can be seen in Figure~\ref{goldengraph}. 
We will describe the function $s(x)$ as an 
explicit series
which can be evaluated in $O(\log(1/\epsilon))$ steps if we want to know it up to 
accuracy $\epsilon$. An analytic function $\sigma$ will be key to determine $s$.  \\

Why is this interesting? Some motivation and references were already given in \cite{KL,KTp}, where we studied the
Birkhoff sum for the function $G$ satisfying $G'=g$. \\

The {\bf Sinai-Ulcigrai model}
$$ g(x) =2 (1-\exp(i x))^{-1} = 1 + i \frac{\sin(x)}{1-\cos(x)}=1+i \cot(\frac{x}{2})  $$ 
is equivalent. In \cite{SU} these authors have shown the existence of a 
limiting distribution of $S_n/n$ if one averages both over the initial point 
$\theta$ and the rotation number $\alpha$. In other words, they have shown that the Birkhoff sum for $\cot$ and 
the map $T(\theta,\alpha)=(\theta+\alpha,\alpha)$ on ${\bf T}^2$ has a limiting distribution in the sense that
the random variables $S_n/n$ converge in law to a nontrivial distribution. 
Where Sinai-Ulcigrai average over all initial conditions and rotation numbers, 
we look at {\bf one} specific orbit and take {\bf one} specific rotation number $\alpha$, 
the golden ratio. We prove that in this particular case, the rescaled random walk 
$S_{q_{2n}}/q_{2n}$ converges.  For the Sinai-Ulcigrai distribution, it is 
not required to go to a subsequence, but one has to integrate over $\theta$ and $\alpha$.  \\

The related aspect of the Sinai-Ulcigrai limiting theorem will be touched upon at the end of the paper, when we look at the
antiderivative $G$ of $g$, the Hilbert transform of the piecewise linear Hecke function $H$ studied
by number theorists. Since Denjoy-Koksma theory \cite{CFS,J} implies that the log Birkhoff sum of $S_n^H/\log(n)$ 
is bounded for almost all $\theta$ and $\alpha$ of constant type, also the Birkhoff sum $S_n^G/\log(n)$ for the 
Hilbert transform $G$ has this property and there is a subsequence $S_{n_k}^G/\log(n_k)$ with a limiting distribution. 
We have explored this limiting distribution for
a single orbit in $S_n^G$ in \cite{KTp}. The result proven in the present paper 
implies that for the golden ratio $\alpha$, the Birkhoff sum $S_{q_n}^{g^{(k}}/q_n^k$ has a 
limiting distribution for derivatives $g^{(k)}, 0 \leq k$, again for a single orbit. \\

In probability theory, the Birkhoff sum under consideration 
deals with examples of ergodic aperiodic sequences of identically distributed
random variables $X_k$ with finite Cauchy principal value expectation.
We can so give explicit expressions for the sum of Cauchy distributed random variables
$$  S_n(\theta) = X_1(\theta)+ \cdots +X_n(\theta) \;  $$
for one point $\theta=0$. 
This is related to Sinai-Ulcigrai, who take as the probability space the two torus and generate the random
variables with the shear $T(\theta,\alpha)=(\theta+\alpha \; {\rm mod} \; 1,\alpha)$ 
and find a limiting distribution when averaging over pairs
$(\theta,\alpha)$ of initial conditions and rotation numbers. Because every natural probability space
is equivalent to the unit interval $\Omega=[0,1]$ with Lebesgue measure $P=dx$, every sequence of 
identically distributed Cauchy random
variables can be obtained as $X_k = \cot(\pi T^k(x))$ with some measure preserving transformation $T: \Omega \to \Omega$. 
Taking an irrational rotation gives stochastic processes with strong correlations between individual random variables.
The "risk" of this process is large because the orbit gets close to the pole, where the "loss" or "gain" can become arbitrary large. 
It complements a trivial "Cauchy central limit theorem" which tells that for any process $X_n = Y_n+Z_n$ with 
identically distributed but not necessarily uncorrelated $L^1$ random variables $Y_k$ 
and IID distributed Cauchy random variables $Z_k$, the sum $S_n/n$ converges in 
law to a Cauchy distributed random variable. This is 'trivial' because $(1/n) \sum_{k=1}^n Y_k$ converges 
to a constant $M$ and because $(1/n) \sum_{k=1}^n Z_k$ all have the same Cauchy distribution. 
The Cauchy distribution shares with the Gaussian distribution the stability of adding independent random variables.
Since the standard deviation is infinite, the Cauchy distribution can model high risk. 
In this paper, we study a limit theorem in the case 
when the Cauchy process has correlations given by a Diophantine rotation. Since the result will be pretty independent
of the function $g$ as long as it has a single pole and the initial point is that pole
it can - like the Sinai-Ulcigrai result - be seen as a 
limit result, where a chaotic Bernoulli process is replaced by a Diophantine irrational rotation.  \\

Theorem~(\ref{main}) allows to estimate $S_n/n$ as $(1/n) \sum_{l=1}^L s(x_l)$, 
where the {\bf Zeckendorff representation} $n = \sum_{k=1}^L a_k q_{k}$ of the integer $n$ is 
linked to the {\bf $\beta$-expansion} $x_l = \sum_{k=1}^l a_k \alpha^{L-k}$. 
Similarly than in the bounded case, we could give an estimate of $S_{10^{100}}$ for example, even so 
there are too many terms to add this up directly. The reason is that we can compute 
$S_n$ for Fibonacci lengths
as $S_n \sim n s(n/q)$, where $q$ is the smallest Fibonacci number larger or equal to $n$.
Given a number $n$ like $n=23'866$, we can represent it as a Zeckendorf sum $23866=17711+4181+1597+377$ of
of Fibonacci numbers $q_{22},q_{19},q_{17},q_{14}$. The next larger Fibonacci number is $46368=q_{24}$. The real number
$x=23866/46368=0.514708$ has the $\beta$ expansion $x=\alpha^2+\alpha^5+\alpha^7+\alpha^{10}$, where $\alpha$
is the golden ratio. We will see that $S_n \sim \sum_{k=0}^{\infty} a_k (-\alpha)^k \sigma(y_k)$
with explicitly known $y_k$. In our example, $y_1=0,y_2=-0.0344396,y_3=-0.026314; y_4=0.00118918$. While the 
actual sum is $-0.187542$, the analytically computed sum obtained by adding up 4 summands is $-0.18763$. 
We have made use of a function $\sigma$ which $n$ independent and analytic and which can be computed at first.
While $S_{10^{100}}/10^{100}$ can not be computed by summation, we could determine the sum analytically. The accuracy
depends on the structure of $n$, on how many Fibonacci numbers are needed to Zeckendorff represent the number.
As the function $\sigma$ is analytic, the Birkhoff sum can so be estimated so fast that it could be called "integrable"
and be diametrally opposite to the IID case, where the value of $S_n/n$ can only be computed in average for 
numbers like $10^{100}$.  \\

For number theory, the Birkhoff sums under consideration
are related to the values of a modular form at a specific point. 
The sums $S_n$ always converge without normalization if $\alpha$ is in the upper half plane
$\{ \alpha \in \C \; | \; {\rm im}(\alpha)>0 \; \}$. The classical {\bf $\theta$ function}
$$   \theta(\alpha) = \sum_{n \in Z} e^{i \pi n^2 \alpha} 
                    = \sum_{n \in Z} w^{n^2} 
                    = 1+2 \sum_{n=1}^{\infty} w^{n^2}  \; $$
with $w=e^{i \pi \alpha}$ is related to the modular form 
$$  R(\alpha)=\theta^2(\alpha)=1+2 \sum_{k=1}^{\infty} g(k \alpha) + g(k \alpha+\pi/2),
                              =1+2 S_{\infty}(\cot,\pi/4) - 2 S_{\infty}(\cot,-\pi/4) $$
for ${\rm im}(\alpha)>0$. In other words, 
the modular form $R$ is up to a constant the sum of two Birkhoff sums described in this article
which converge for $\alpha$ in the upper half plane.  \\

The sums also relate to the {\bf theory of partitions} in additive number theory because $g(x)=G'(x)$ with 
$G=\log(1-z)$ and $z=e^{2\pi i x}$. The Birkhoff sum for $G$ is the logarithm of a product 
$$  P_n=\prod_{k=1}^n (1-w^k)  \; , $$
where $w=e^{2\pi i \alpha}$. The function $-G$ is related to the generating functions for the partition function $p(n)$ because 
$P_n^{-1}= \sum_{k=1}^n p(k) w^k + O(w^{n+1})$.

To illustrate the connection between number theory and dynamical systems,
lets prove a well known statement in the theory of {\bf integer partitions} $p(n)$ (which tells in 
how many ways we can write $n$ as a sum of positive integers) with 
the help of Birkhoff sums in dynamical systems theory: \\
\begin{center}
{\it The partition function $p(n)$ satisfies $\limsup_n |p(n)|^{1/n} =1$.}
\end{center}
\begin{proof}
Euler's {\bf pentagonal number theorem} (i.e. \cite{SteinshakarchiComplex})
shows that the function $Q(z) = \prod_{k=1}^{\infty} (1-z^k) = \sum_n a_n z^n$ 
which is related to the partition numbers by $P(z)=\sum_{n=1}^{\infty} p(n) z^n = \prod_{k=1}^{\infty} (1-z^k)^{-1} = 1/Q(z)$
has Hadamard gaps. By Hadamard's {\bf lacunary series theorem}, the Taylor series has a natural boundary as a circle.
This implies that establishing boundedness on one point $|z|=r$ implies that the radius of convergence is larger 
than $r$. Define $g(\alpha) = (1-r e^{2\pi i\alpha})^{-1}$ so that $f(\alpha) = \log(P(z)) = - \sum_{k=1}^{\infty} g(k \alpha)$ 
For $r=|z|<1$, and $z=\exp(2\pi i \alpha)$ with Diophantine
$\alpha$, the theorem of Gottschalk-Hedlund \cite{KH} shows that the function $f$ is a coboundary 
$f(x) = F(x+\alpha)-F(x)$ implying that the sum is bounded. Having so proven that the Taylor series has no singularity on
$|z|<1$,the partition function $p(n)$ satisfies $\limsup_n |p(n)|^{1/n} \geq 1$. 
The other inequality $\limsup_n |p(n)|^{1/n} \leq 1$ is trivial because
$p(n) \geq 1$ implies $\sum_{n=1}^{\infty} p(n) z^n \geq \sum_{n=1}^{\infty} z^n$.
\end{proof}

While this growth-statement for the partition numbers $p(n)$ is not difficult to prove directly \cite{Newman},
the just given new argument illustrates how complex analysis, dynamical systems theory and Diophantine notions can 
play together. \\

\begin{figure}
\scalebox{0.9}{\includegraphics{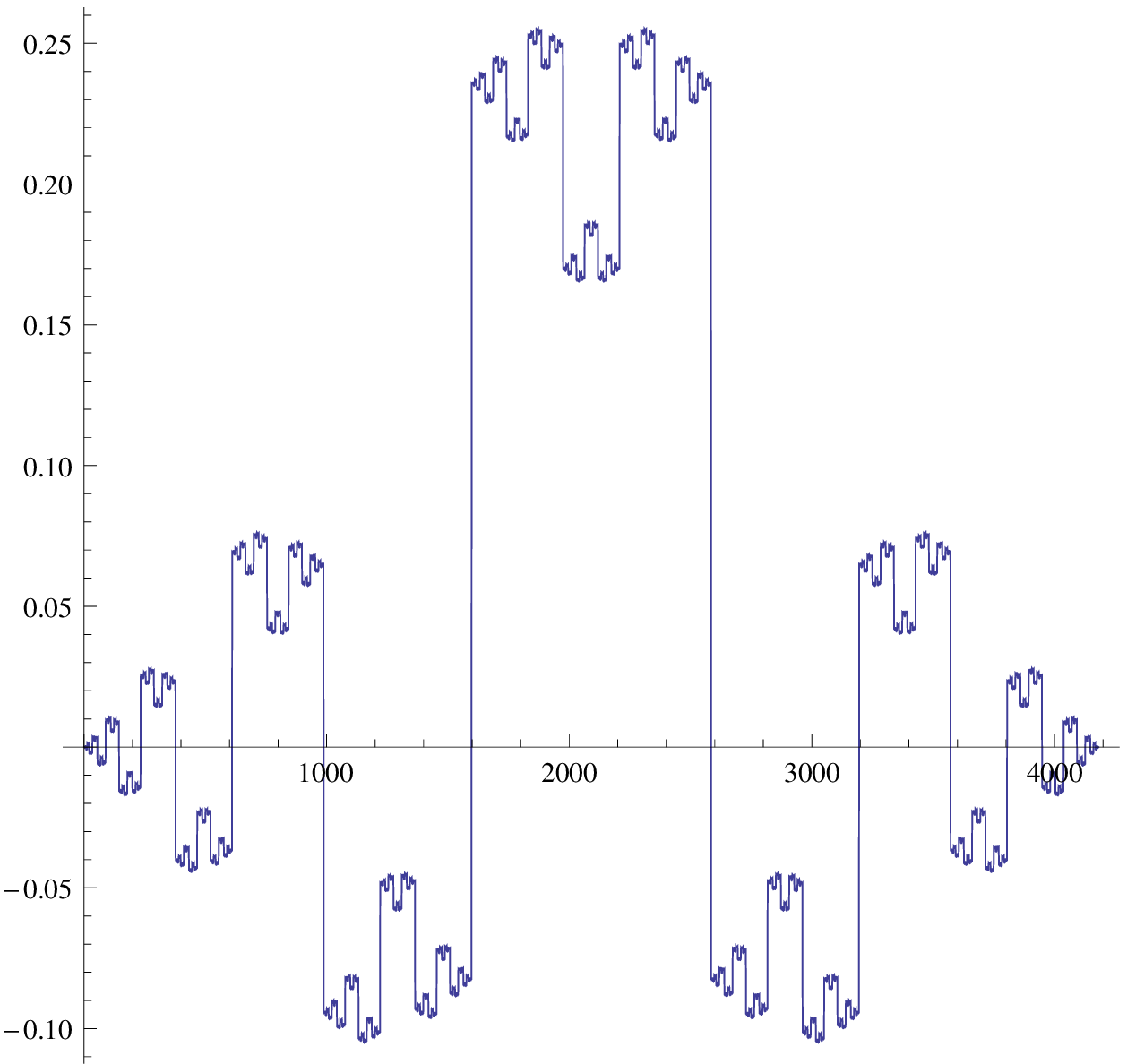}} \\
\caption{
The {\bf rational golden graph} is the Birkhoff limiting function
$$  t(x) =   \lim_{n \to \infty} \frac{S_{[q_{2n} x]}}{q_{2n}}(T_n) $$
in which the translation $T_n(x)=x+p_{2n}/q_{2n}$ is periodic. 
It is different from the Birkhoff limiting function $s$ in which $T(x)=x+\alpha$ is $n$ independent. 
The function $t$ is not affine self similar 
($a+b(-t(\alpha x)-t(\alpha(1-x))$ only comes close).
But the function $t$ can be useful to answer questions posed in \cite{KTp}. 
}
\label{cotsumrational}
\end{figure}

The choice of the core function $g(x)=\cot(\pi x)$ for the Birkhoff sum is distinguished 
because the periodic function $g$ has {\bf as a distribution} the constant Fourier expansion
\begin{equation}
\label{constantfourier}
g(x) = \cot(\pi x) = 2 \sum_{k=1}^{\infty} \sin(2\pi k x)  \; 
\end{equation}
as can be seen by taking the real part of 
$$ \cot(x) = i \frac{e^{ix}+e^{-ix}}{e^{ix}-e^{-ix}} = -i \frac{1+e^{2i x}}{1-e^{2ix}} = -i (1+2 \sum_{k=1}^{\infty} e^{2 i kx}) \; . $$
How special the cot function is can also be seen with the {\bf Euler's cot formula}
\begin{equation}
g(x) = \frac{1}{\pi} \sum_{k \in Z} \frac{1}{x-k}  \; .
\end{equation}
Appearing also as the Hilbert kernel it is the most natural periodic function with a single pole. \\

In \cite{KL,KTp}, we studied the Birkhoff sum of the antiderivative
of $g(x) = \cot(\pi x)$, which is
\begin{eqnarray*}
    G(x) &=& \log(2-2 \cos(2\pi x)) \frac{1}{2\pi}  \\
         &=& \log(2 \sin(\pi x)) \frac{1}{\pi}  \\
         &=& \log|1-e^{2\pi i x}| \frac{1}{\pi}  \; . 
\end{eqnarray*}
Euler's formula for $g(x)=G'(x)$ can also be deduced by logarithmic differentiation 
of the {\bf Euler's product formula} for the sinc-function
$$  {\rm sinc}(\pi x) = \frac{\sin(\pi x)}{x \pi} = \prod_{k \neq 0} \frac{x-k}{k} $$ 
because $2 \sin(\pi x) = \exp(\pi G(x))$. \\

The current paper will lead to explanations to some of the statements found in the experiments done in \cite{KTp}.
Wee do not give much details here yet but say that one connection is the fact that
the difference of the Birkhoff sums of the antiderivative $G(x)$ 
of $x$ for successive continued fractions $p_n/q_n, p_{n+1}/q_{n+1}$ can be expressed by Rolle's theorem as
a Birkhoff sum of $g'$ but with a different initial condition. While the golden graph of $g$ is selfsimilar, the
limiting distributions seen in \cite{KTp} are only almost self similar but the difference can be expressed by a 
Birkhoff sum for $g'$ which is monotone. \\

While the rotation number $\alpha$ is special because it has a constant continued fraction expansion,
the story could be generalized a bit. For quadratic irrationals, there is a periodic 
attractor for the limit (see \cite{KTp}). \\

For $L^1$ functions $f$, by the {\bf Birkhoff's ergodic theorem},
the limiting function is for almost all $\theta$ given by
$$ s(x) 
        = \lim_{n \to \infty} \frac{1}{n} \sum_{k=1}^{[x n]} f(T^k \theta) 
        = \lim_{n \to \infty} \frac{x}{m} \sum_{k=1}^{m} f(T^k \theta) 
        = M x \; , $$
where $M = \int f(\theta) \; d\theta$. We are not aware of any other case of a stochastic process,
where the Birkhoff average $s_n(x) = S_{[n x]}/n$ converges to a nonlinear function on $[0,1]$ along an orbit.
The $\cot$ function can be modified by adding any continuous periodic function $F$:

\begin{coro}[A single orbit ergodic theorem for the golden rotation]
For any function $g$ on $T^1=[0,1)=R/Z$ which is continuous away from $\theta=0$ 
and has a single pole at $\theta=0$, the Birkhoff limiting function at $\theta=0$ 
$$  s(x) = \lim_{n \to \infty} \frac{S_{[q_{2n} x]}}{q_{2n}} $$ 
is the sum of a multiple of the golden graph $c \cdot s(x)$ and a linear function $Mx$.
\end{coro}

\begin{proof}
Because $f$ has a single pole at $0$, we can write $f(x) = c \cot(\pi x) + F(x)$, where $F(x)$ is continuous and 
especially in $L^1$.  The Birkhoff limiting function for $f$ is the sum of the golden graph $c s(x)$ 
and the Birkhoff limiting function of $F$, which is $M x$ with $M=\int_0^1 F(x) \; dx$.
\end{proof}

\begin{figure}
\scalebox{0.9}{\includegraphics{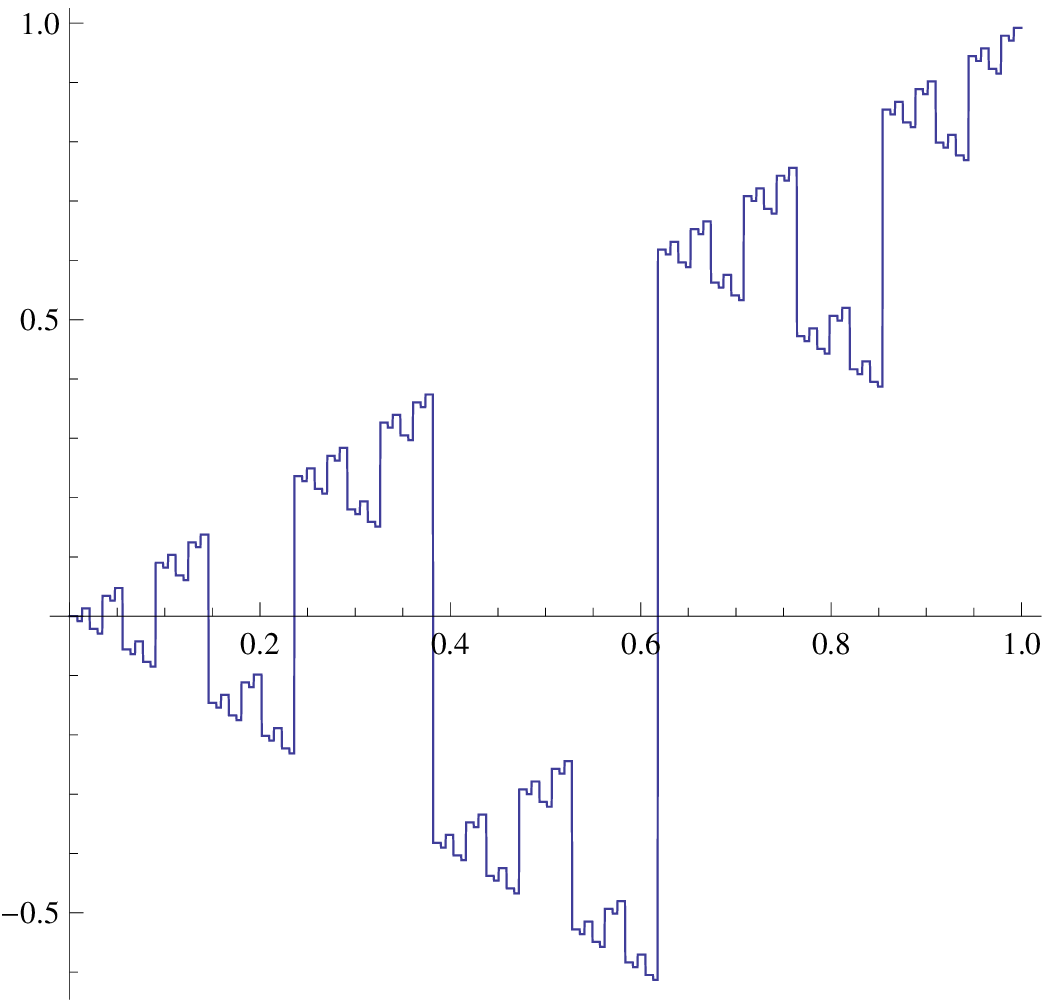}} \\
\caption{
A {\bf caricature for the golden graph} is the function
$$ o(x) = \sum_{k=1}^{\infty} (-1)^{k-1} a_k \alpha^k \; , $$
where $\alpha=(\sqrt{5}-1)/2$ and $a_k$ is the $\beta$-expansion of $x$.
The function $o$ has the same symmetry $o(\alpha^n x) = (-\alpha)^n o(x)$
as the golden graph. The golden graph is a modification of this graph since we will show
$$ s(x) = \sum_{k=1}^{\infty} a_k (-\alpha)^k \sigma(y_k) \; , $$
where $y_0=0$, $y_k= \lim_{n \to \infty} [|x_k q_n| \alpha]q_n$ with 
$x_k = \sum_{i=1}^k a_i \alpha^i$. 
}
\label{cotsum}
\end{figure}

We will be able to compute the function $s(x)$ explicitly from a function
$$  \sigma(y) = \lim_{n \to \infty} \frac{1}{q_n} \sum_{k=1}^{q_n-1} g(\frac{y}{q_n} + k \alpha) \; $$ 
which we will show to be analytic in $y$. Its Taylor expansion starts with 
$$  \sigma(y) = 0.258 + -1.24724 y + 0.6736 y^2 -1.38082 y^3  + ...  \; .  $$

First, we look at 
the analogue function for {\bf rational} $\alpha$, where the function is 
$$  \tau(y)  =  \lim_{n \to \infty} \frac{1}{q_n} \sum_{k=1}^{q_n-1} g(\frac{y}{q_n} + k \frac{p_n}{q_n}) \; . $$
The function $\tau$ is odd and will be given explicitely: 
$$ \tau(y)   =   \frac{\pi}{3} y - \frac{\pi^3}{45} y^3 - \frac{2\pi^5}{945} y^5 + \dots  \; . $$

From $\sigma$ we can then get $s$. \\

An indication why things work nicely is that the $\cot$ function and their derivatives 
have structure. Hardy and Littlewood have made use of this early on. 
The selfsimilar nature of the Birkhoff sum appears in the following statement which we have
used already in integrated form in \cite{KL} and which will be used in an equivalent way below in 
Lemma~\ref{explicitrational}. 

\begin{propo}[Birkhoff renormalization]
For any periodic translation $T_{p/q}(y)=y+p/q$ on the circle $R/Z$ with ${\rm gcd}(p,q)=1$,
the function $g(y) = \cot(\pi y)$ is a fixed point $B_{p/q}(g)=g$ of the linear
{\bf Birkhoff renormalization operator} 
$$  B_{p/q}(g)(y)=\frac{1}{q} \sum_{k=0}^{q-1} g(T_{p/q}^k (\frac{y}{q})) \; . $$
\label{cotrenormalization}
\end{propo}

{\bf Remark.} \\
We believe that any periodic function which is a fixed point of all the Markov 
operators $B_{p/q}$ must be of the form $f(y) = a+b \cot(\pi y)$ with some 
constants $a,b$. But we have not yet been able to prove such a fixed point result. 

\begin{figure}
\scalebox{0.9}{\includegraphics{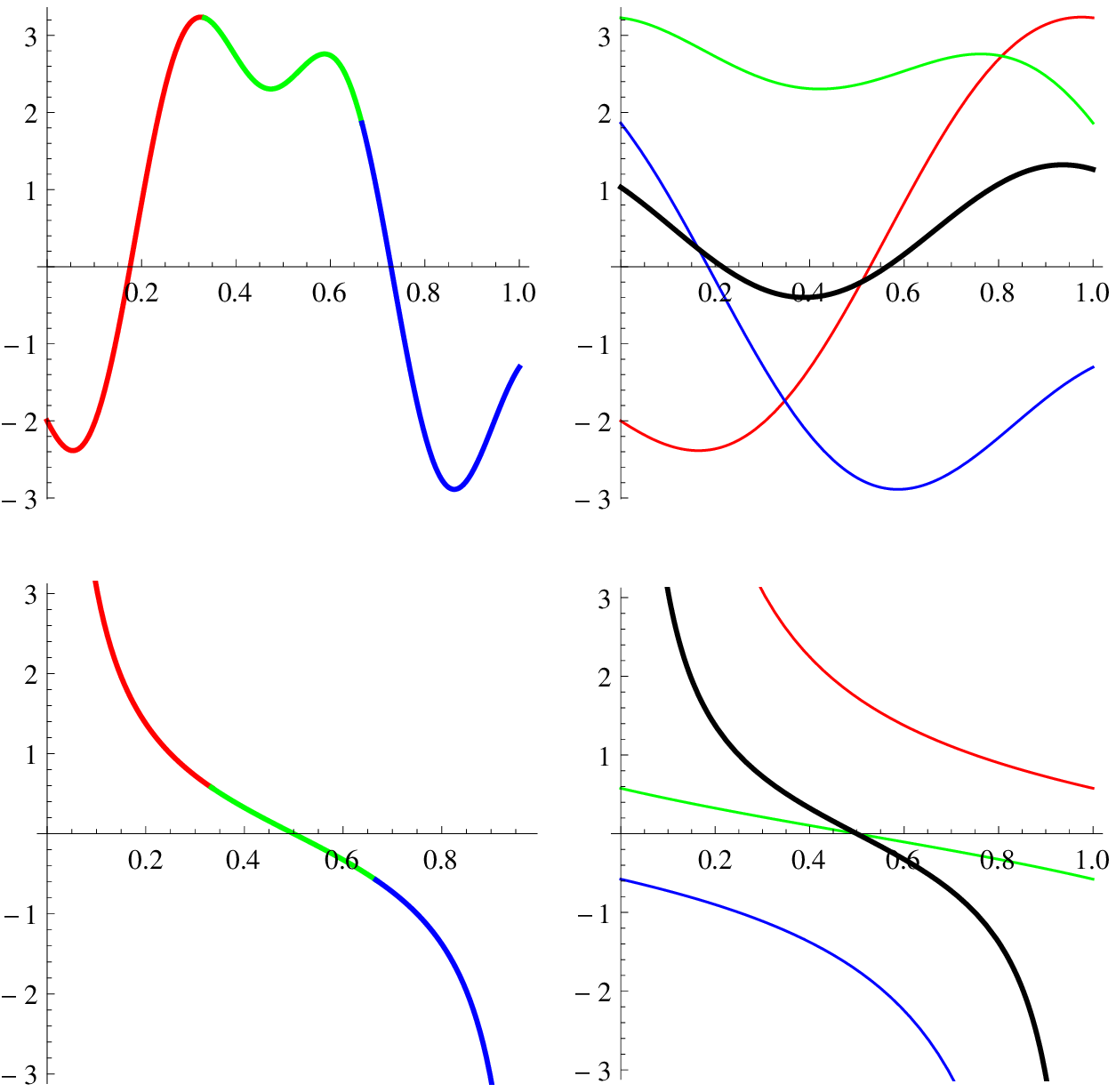}} \\
\caption{
The Markov operator $B_{1/3}$ is illustrated in the figure as it is applied to a general function, where
it produces a new different function. The two figures  below show the operator applied to the $\cot(\pi x)$ function
which is a fixed point: 
$$  \frac{1}{3} [ \cot(\pi \frac{y}{3}) + \cot(\pi \frac{y}{3} + \pi \frac{1}{3}) 
                                        + \cot(\pi \frac{y}{3} + \pi \frac{2}{3}) ] = \cot(\pi y) \; . $$
We see the three graphs, obtained by taking stretched pieces of the cot graph and then 
averaging them to get back the function. Each operator $B_{p/q}$ is a Perron-Frobenius operator 
for an interval exchange transformation.
}
\label{cotsum}
\end{figure}

\begin{proof}
The operators $B_{p/q}$ are Markov operators because they satisfy $B 1=1$ and  $||Bf||_1 = ||f||_1$
for positive $L^1$ functions. Use the {\bf cyclotomic formula}
$$ \prod_{k=1}^{n-1} (1-e^{2\pi i k/n} z) = \frac{z^n-1}{z-1} \;  $$
and plugging in $z=\exp(2\pi i y)$, then differentiate with respect to $y$ to get
$$ \frac{1}{q} \sum_{k=0}^{q-1} \cot(\pi \frac{(k+y)}{q}) = \cot(\pi y)   \; .  $$
This is equivalent because $\{ p k \; {\rm mod} \; q \; \} = \{ 0, \dots , q-1 \; \}$ follows from
${\rm gcd}(p,q)=1$. We have so shown that $g=a+b \cot(\pi x)$ is a fixed point of $B_{p/q}$
for all $p/q$. 
\end{proof}

{\bf Remarks.} \\
{\bf 1)} Unlike for other parts of this article which often were found experimentally, 
the fixed point property was obtained first in a deductive way from the cyclotomic identity.  \\
{\bf 2)} Because related identities appear in formula collections,
we must assume Proposition~\ref{cotrenormalization} to be "known". Still, it is a 
remarkable gem and probably not "well known". Here are some related identities
$(1/q^2) \sum_{k=0}^{q-1} \csc^2(\pi (k+y)/q) = \csc^2(y)$ or
$(1/q^2) \sum_{k=0}^{q-1} \cot^2(\pi (k+y)/q) = \cot^2(y) + 1-1/q$ for $\cot^2$
which are equivalent by differentiation and using the trig identity $\cot^2(x)+1=\csc^2(x)$. 
We get so functions which are fixed points of the higher order Birkhoff operators
$$  B_{p/q}^{(l)}(g)(y)=\frac{1}{q^l} \sum_{k=0}^{q-1} g(T_{p/q}^k (\frac{y}{q}))=g(y) \; . $$
{\bf 3)} In the simplest case, the renormalization picture provides for $p/q=1/2$ a double angle formula 
for the $\cot$ function: for $x=\pi y/2$ it is 
$$ \cot(x) + \cot(x+\frac{\pi}{2}) = 2 \cot(2x)  \; . $$
For other $p/q$, we get {\bf "de Moivre type" trigonometric identities} for the $\cot$ or $\tan$
functions analogues to the magic which de Moivre has provided for $\cos$ and $\sin$. For example
$$ \cot(x)+\cot(x+\frac{\pi}{5})+\cot(x+\frac{2\pi}{5})+\cot(x+\frac{3\pi}{5})+\cot(x+\frac{4\pi}{5})=5 \cot(5x) \; . $$

Before we launch into the details of the Birkhoff sum puzzle, here is a schematic overview over
the proof. We start with an explicit Birkhoff sum in the rational case by using an elementary
cyclotomic formula leading to the function $\tau(y)$. Interpolating this with Fibonacci
steps allows us to get a related analytic function $\sigma(y)$ for which we can compute the Taylor
coefficients with arbitrary precision. This in turn 
provides us with a function $s(x,y)$ via a $\beta$ expansion of $x$. The slice $s(x)=s(x,0)$ 
is then the golden graph.  This golden graph in turn allows us via Taylor expansion and Abel summation
to compute analytic "golden graphs" in the rational case $t(x,y)$. This function $t$ is not selfsimilar 
but it will link us to the original observations done in \cite{KTp}. 
It is interesting to note that the direct way in the rational case from $\tau$ to $t(x,y)$ is blocked
because in the rational case, there is no self similarity. We need to make a detour over the 
golden ratio and the golden graph to understand the rational case!

\begin{center}
\begin{pspicture}(0,0)(18,4)
\rput(3,4){\fb{$\tau(y)=t(1^-,y)$}}   
\rput(3,1){\fb{$\sigma(y)=s(1^-,y)$}}
\rput(7,4){\fb{$t(x,y)$}}     
\rput(7,1){\fb{$s(x,y)$}}    
\rput(11,1){\gfb{$s(x)=s(x,0)$}}    
\rput(11,4){\fb{$t(x)=t(x,0)$}}    
\psline[linearc=.25]{->}(8,4)(9,4)
\psline[linearc=.25]{->}(5,1)(6,1)
\psline[linearc=.25]{->}(8,1)(9,1)
\psline[linearc=.25]{->}(3,3)(3,2)
\psline[linearc=.25]{->}(7,2)(7,3)
\rput(3,0.2){(Fig 7)}
\rput(7.3,0.2){(Fig 6)}
\rput(10.7,3.0){(Fig 2)}
\rput(10.9,0.2){Golden Graph (Fig 1)}
\rput(5.9,0.2){$\beta$-exp.}
\rput(8,2.5){Taylor}
\rput(1.9,2.5){Steps}
\end{pspicture}
\end{center}

{\bf Acknowledgements:}
Thanks to John Lesieutre and Folkert Tangerman for previous 
joint work [KL] and [KTp] without which this paper would not have been possible. 

\section{Constructing the function $s$}

Every real number $x \in [0,1]$ has a {\bf $\beta$-expansion} $\sum_{n=1}^{\infty} a_n \alpha^n$ with $a_i \in \{0,1 \}$. 
If $\alpha=(\sqrt{5}-1)/2$, then the possible sequences $\{a_n \; \}_{n}$ are the ones for which no 
two consecutive $a_i$ are equal to $1$ (see for example \cite{DK}). 
This is related to the {\bf Zeckendorf representation} $n = \sum_{n=1}^{\infty} a_n q_n$
of integers as a sum of Fibonacci numbers, where only finitely many $a_n$ are nonzero. Both in the $\beta$
expansion of a real number in $[0,1]$ as well as the Zeckendorf representation of a positive integer,
two successive $1$'s are forbidden in the subshift of finite type of all $\{0,1 \; \}$ sequences 
because $\alpha^k+\alpha^{k+1}=\alpha^{k-1}$ and so $q_{n-1}+q_n=q_{n+1}$. 
Every point $x \in [0,1)$ has a $\beta$ expansion. Examples of some $\beta$-rational points are
$x=\alpha=0.1,\alpha^2=0.01$, $\alpha+\alpha^3 = 0.101_{\alpha}$ or 
$\alpha+\alpha^6+\alpha^8 = 0.10000101_{\alpha}$. \\

\begin{figure}
\scalebox{0.9}{\includegraphics{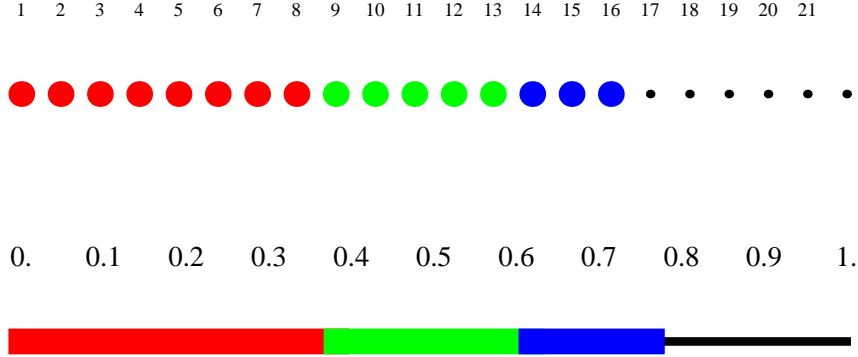}} \\
\caption{
The Zeckendorff representation of a number is the discrete analogue or dual to the 
$\beta$ expansion of a real number. The picture illustrates the Zeckendorff 
representation of $16=3+5+8 = q_2+q_3+q_4$ as a sum of Fibonacci numbers as well as the $\beta$ expansion of 
$0.763932... = \alpha^2+\alpha^3+\alpha^4$ as a sum of powers of the golden ratio.
}
\label{cotsum}
\end{figure}

The semi-continuity of $s$ is reflected by the fact that  
for every $k>0$, the sum $S_{k+q_n}/q_n$ has the same limit. For $k<0$, the limit is different.
To understand the limit, it is useful to start at different points and 
to look more generally at the Birkhoff sum
$S_n(\theta) = \sum_{k=1}^{n} g(\theta + k \alpha)$. We choose $\theta=y/q_{2n}$ 
and define the {\bf $s$-function}
$$  s(x,y) =   \lim_{n \to \infty} S_{[q_{2n} x]}(y/q_{2n})/q_{2n} 
           = - \lim_{n \to \infty} S_{[q_{2n+1} x]}(y/q_{2n+1})/q_{2n+1} \; .  $$

\begin{figure}
\scalebox{0.9}{\includegraphics{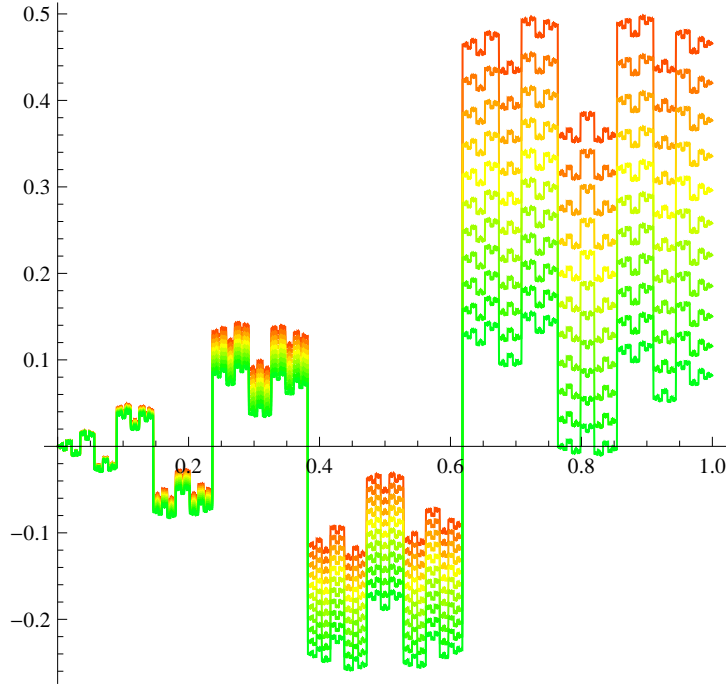}} \\
\caption{
The generalized golden graphs $x \to s(x,y)$ are seen here for different initial points $y$. 
The figure shows the graphs for $y \in [-0.15,0.15]$. We prove in this paper that
the graphs deform continuously
in $y$ for small $|y|$ (one could actually verify it for $|y| <1/\sqrt{5}$). With the initial conditions 
$y_n=y/q_n$ with small $|y|$, the limit
$S_{[x q_n]}(y_n) \to s(x,y)$ still exists along even and odd subsequences. 
}
\label{cotsum}
\end{figure}

We show now that the function 
$$  \sigma(y) = s(1^-,y) = \lim_{n \to \infty} S_{[q_{2n+1}]}(\frac{y}{2q_n+1})/q_{2n+1}    $$  
determines $s(x,y)$ and that the graph of $s$ is selfsimilar. We get

\begin{propo}[The golden graph in terms of $\sigma$]
\label{sigmatos}
Assume we know that $\sigma(y)=s(1,y)$ exists. 
If $x=\sum_{j=1}^\infty a_j \alpha^j$, then
$$ s(x,y) = \sum_{k=0}^{\infty} a_k (-\alpha)^k \sigma(y_k) \; , $$
where $y_0=y$, $y_k=y + \lim_{n \to \infty} [|x_k q_n| \alpha]q$ depend on $x$ via
$x_k=\sum_{j=1}^k a_j \alpha^j$.
We have for every $n$ the symmetry $s(\alpha^n x,y) = (-1)^n \alpha^n s(1,y)$.
\end{propo}

\begin{figure}
\scalebox{0.6}{\includegraphics{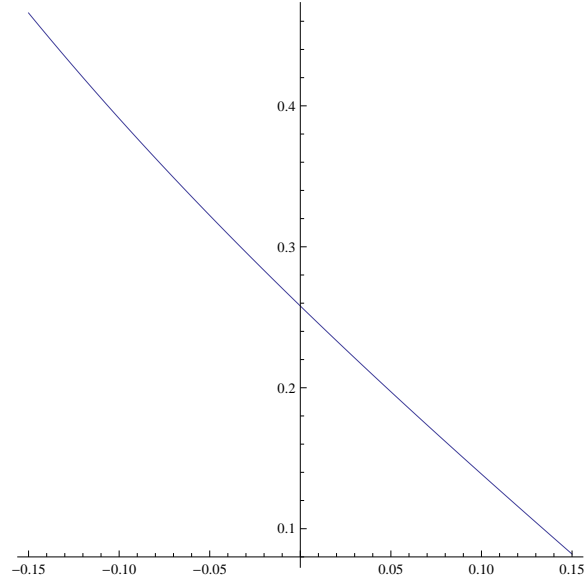}} \\
\caption{
The analytic function $\sigma(y)$ determines the function $s(x,y)$ by a $\beta$ 
expansion as shown in Proposition~\ref{sigmatos}.}
\label{sigmam}
\end{figure}

\begin{proof}
Given $x = \sum_{i=1}^{\infty} a_i \alpha^i$ define the real numbers $x_k = \sum_{i=1}^k a_i \alpha^i$. 
The selfsimilarity claim follows from $s(\alpha x,y) = -\alpha s(1,y)$ which is
a consequence of the fact that for the golden ratio, 
$$  \alpha q_n-q_{n+1} = \frac{\pm 1}{q_n} \; . $$

Use that for $x=\alpha^k$, the selfsimilarlity implies that 
$s(\alpha^k,y) = (-\alpha)^k \sigma(y (-\alpha)^k)$. Use the notation $[x] = x \; {\rm mod} \; 1$. 
Define $y_0 = y$ and $y_k=y + \lim_{n \to \infty} |x_k q_n| \cdot \alpha \cdot q_n$. Now,
$$ s(x,y) = \sum_{k=1}^{\infty} s(\alpha^k,y_k) = \sum_{k=1}^{\infty} (-\alpha)^k \sigma((-\alpha)^k y_k) \; . $$ 
\end{proof}

The main task is therefore to construct the function $\sigma$. 
To do so, we have first to explore the rational case, where we sum up $q_n-1$ terms
with the rotation number $p_n/q_n$ instead of with $\alpha$. The analogue to 
the $\sigma$ function is called $\tau$. It can be computed explicitly for $g(x) = \cot(\pi x)$ and defined as
$$ \tau_{p/q}(y) = \frac{1}{q} \sum_{k=1}^{q-1} g(\frac{y}{q}+\frac{k}{q}) \; . $$
This expression can be simplified and also explicitly be given in the limit $q \to \infty$. 

\begin{lemma}[The function $\tau$] 
\label{explicitrational}
If $\{ k p \; \}_{k=1}^{q-1} = \{ 1, \dots, q-1 \; \}$ then 
$$ \tau_{p/q}(y) = g(y) - \frac{g(\frac{y}{q})}{q} \; . $$
In the limit $q \to \infty$,
$$  \tau(y) = g(y) - \frac{1}{\pi y} \; . $$
\end{lemma}
\begin{proof}
The {\bf cyclotomic formula}
$$ \prod_{k=1}^{n-1} (1-e^{2\pi i k/n} z) = \frac{z^n-1}{z-1} \;  $$
holds for all complex numbers $z \neq 1$. Taking logarithms and looking at the real part,
we obtain for $z = e^{2\pi i y}$ and real nonzero $y$ the formula
$$ \sum_{k=1}^{n-1} \log|2-2 \cos(2\pi \frac{k}{n} 
     + 2\pi \frac{y}{n})| = 2 \log| \frac{e^{2\pi i y}-1}{e^{2\pi i y/n}-1}|
                          = G(y)  - 2 \log|e^{2\pi i \frac{y}{n}}-1| $$
so that 
$$ \sum_{k=1}^{n-1} G(\frac{y}{n}+\frac{k}{n}) = G(y) - G(\frac{y}{n})      \; . $$
Now differentiate with respect to $y$ to get
$$ \frac{1}{n} \sum_{k=1}^{n-1} g(\frac{y}{n}+\frac{k}{n}) = g(y) - g(\frac{y}{n})/n  \; . $$
Now use $g(y/q) \sim q/(y \pi)$ to get the statement in the lemma.
\end{proof}

{\bf Remark.} Repeated differentiation gives more identities
$$ \lim_{n \to \infty} \frac{1}{n^{l+1}} \sum_{k=1}^{n-1} g^{(l)}(\frac{y}{n}+\frac{k}{n}) 
   \to g^{(l)}(y) - g^{(l)}(\frac{y}{n})/n^l  \; . $$

For later, we also will need the rational analogue of the $s$ function which we call the {\bf t-function}
$$ t(x,y) = \lim_{q \to \infty} \frac{1}{q} \sum_{k=1}^{[q x]} g(\frac{y}{q} + k \frac{p}{q}) \; . $$
We see its graph in Figure~(\ref{cotsumrational}) for $y=0$.  We know that $t(1,y) = \tau(y)$
is analytic, as we have given it explicitly. 
The function $t$ is only close to self similar but does not quite fit an affine similarity
$t(\alpha x,r y = t(x,y)a + b$. Therefore, we can not get $t$ from $\tau$ as we did get $s$ from $\sigma$. 
In the rational case, we have an explicit limiting function for $\tau$ $x=1^-$ but no 
self similarity. In the irrational case, we have a self similar $s$ function but no explicit expressions
for $\sigma$. 

\begin{coro}[The $t$-function]
If $s(x,y)$ is continuous to the right and analytic in $y$ 
then the function $t(x,y)$ is continuous to the right in $x$ and analytic in $y$.
\end{coro}

\begin{proof}
We can get $t(x,y)$ from the function $s(x,y)$ by making a Taylor expansion in $\alpha$ centered at the golden ratio
and evaluate it at a continued fraction approximations $p/q$: 
$$ t_{p/q}(x,y) = s(x,y) + s'(x,y) (\frac{p}{q}-\alpha)/1! + s''(x,y) (\frac{p}{q}-\alpha)^2/2! + \cdots  \; . $$
and then take the limit $p/q \to \alpha$. 

We have $|\frac{p}{q}-\alpha| \leq \frac{1}{\sqrt{5} q^2}$ and 
$$ A_m = \frac{d^m}{dy^m}  g(\frac{y}{q} + k \alpha)| = \frac{k^m}{q^m} g^{(m)}(\frac{y}{q} + k \alpha)  $$
We need to show that the Taylor series for $t_{p/q}(x,y)$ converges. For fixed $q$, it can be written as 
a sum $\sum_{m} A_m/m!$ which has terms which satisfy for large $q$ (neglecting $1/\sqrt{5}<1$)
$$ A_m \sim  \frac{1}{q^{m}} \sum_{k=1}^{[q x]} \frac{k^m}{q^m} g^{(m)}(\frac{y}{q} + k \alpha)) =\frac{1}{q^{m}} \sum_{k=1}^{[q x]} h_k g_k \; .   $$
With $G_l = \sum_{j=1}^l g_j$ we can write the sum using Abel summation (discrete integration by parts) as 
$$          \sum_{k=1}^{[q x]} h_k g_k  = G_{[qx]} h_{[qx]} - \sum_{k=1}^{[qx]} G_k (h_k-h_{k-1})  \; .  $$
Now
$$ \frac{1}{q^{m}} \sum_{k=1}^{[qx]} G_k (h_k-h_{k-1})  
 = \frac{1}{q^{m}} \sum_{k=1}^{[qx]} G_k (\frac{k^m}{q^m}-\frac{k^{m-1}}{q^{m}}) = O(\frac{1}{q}) $$
so that the second part after Abel summation goes to zero. With $(1/q^m) G_{[q x]} \to s^{(m)}(x,y)$ and 
$\frac{h_{[q x]} }{q^{m}} = (q x)^m/q^m = x^m$ we finally see $A_m \sim s^{(m)}(x,y) x^m$ for $q \to \infty$. 
\end{proof}

The following proposition gives $\sigma$ from the analytic function $\tau$: 

\begin{propo}[From $\tau$ to $\sigma$]
The function $\sigma(y)$ is analytic in $y$ for small $|y|$ 
with Taylor expansion $\sigma(y) = \sum_{l=0}^{\infty} a_l y^l/l!$, where 
$$  a_l = \lim_{n \to \infty} \frac{1}{q^{l+1}} \sum_{k=1}^{q-1} k g^{(l)}(\frac{y}{q}+k \alpha)  $$ 
and $q=q_n$ is the $n$'th Fibonacci number.
\end{propo}

\begin{proof} 
We will compare the functions $\sigma$ to $\tau$ which are given as the limits
\begin{eqnarray*}
 \sigma_n(y) &=&  \frac{1}{q} \sum_{k=1}^{q-1} g(\frac{y}{q}+k \alpha)  \\
 \tau_n(y)   &=&  \frac{1}{q} \sum_{k=1}^{q-1} g(\frac{y}{q}+k \frac{p}{q})   \; 
\end{eqnarray*}
and show that the difference is analytic. To do so, we write the difference as a sum
of differences of higher order rational tau functions, then linearize
$$ \sigma_n(y) - \tau_n(y) 
            = \sum_{l=0}^{\infty} \tau_{n+2l+2}(y)-\tau_{n+2l}(y) 
            = \sum_{l=0}^{\infty} \tau_{n,2l+1}'(y)(\frac{P_{l+2}}{Q_{l+2}} - \frac{P_l}{Q_l})  + R_n\; , $$
with 
$$ \tau_{n,l}'(y)  =  \frac{1}{q} \sum_{k=1}^{q-1} \frac{k}{Q_l} g'(\frac{y}{q}+k \frac{P_l}{Q_l})  \; , $$
where $q=q_n, p=p_n, Q_l = q_{n+l},P_l=p_{n+l}$ for the sake of notational simplicity. The Taylor rest term
$$ R_n  = \sum_{l=0}^{\infty} \tau_{n,2l+1}^{(2)}(y)(\frac{P_{l+2}}{Q_{l+2}} - \frac{P_l}{Q_l})^2/2 $$
disappears in the limit $n \to \infty$ because $(\frac{P_{l+2}}{Q_{l+2}} - \frac{P_l}{Q_l} )
=( \frac{1}{Q_l Q_{l+2}})^2 \sim \alpha^2 \frac{1}{Q_l^4} \sim \alpha^{2+4l} \frac{1}{q_n^4}$
so that remark after Lemma~\ref{explicitrational} shows $R_n \to 0$. 
Because $\tau_{n,l}'(y)$ is a sum of negative terms, the
boundedness of $\tau_{n,l}'(y)$ follows from the Cauchy-Schwarz inequality.
We know that the limit $n \to \infty$ of each $\tau_{n,l}'(y)$ exists as an analytic function 
and that it is bounded above by $C \alpha^l$. Therefore also the sum converges and is an analytic function.
It is the difference between $\sigma$ and $\tau$ we were looking for. 
\end{proof}

{\bf Remarks.}  \\
{\bf 1)} It is in this proof that we see the reason why we have to go along the subsequence $q_{2n}$ to get the 
golden graph and not $q_n$. The reason for the sign flip is: 
$$ \frac{p_{n+2}}{q_{n+2}} - \frac{p_n}{q_n} = \frac{+1}{q_n q_{n+2}},  n \; {\rm odd}   $$
$$ \frac{p_{n+2}}{q_{n+2}} - \frac{p_n}{q_n} = \frac{-1}{q_n q_{n+2}},  n \; {\rm even}  $$
The $\tau$ function has an analytic formula for any $k$, we do not even need a Fibonacci number. 
It is the difference between $\tau$ and $\sigma$ which flips the sign. \\
{\bf 2)} A Taylor expansion in $\alpha$ does not work when going from $\tau$ to $\sigma$. (We only make Taylor expansions
with respect to $y$.) Abel summation indicates that this would require us to know Cesaro sums of the function $t(x,y)$.
The Taylor expansion will work however in order to get $t$ from $s$.
We have written therefore $\sigma-\tau$ as a sum $\sum_{l} f_{l,l+2}$ of smaller differences and estimated
this as a Birkhoff sum with rotation number $p_{n+1}/q_{n+1}$. We can estimate each entry to be 
smaller than $1/q_l$ and so get a Cauchy sequence. There is no need to use Cesaro convergence
because the Birkhoff sums for $g'$ contain only negative entries and establishing boundedness is enough.  \\

While we only need the first approximation, lets state 
a lemma generalizes approximations like $\sqrt{5} (\alpha-p/q)  = (-1)^{n+1}/q^2 + O(q^{-4})$ 
and $\sqrt{5} (\alpha-p/q) - (-1)^{n+1}/q^2 = 1/(5 q^4) + O(q^{-6})$ etc.
It is an exact formula estimating the error from the golden ratio $\alpha$ to any of its continued fraction 
approximations $p_n/q_n$:

\begin{lemma}[Catalan gold]
For any partial fraction $p_n/q_n$ of $\alpha=(\sqrt{5}-1)/2$ we have
$$ \sqrt{5} (\alpha- \frac{p_n}{q_n}) = \sum_{k=0}^{\infty} 
   \frac{(-1)^{(n+1)(k+1)} c_k }{5^k} \frac{1}{q_n^{2k+2}} = (-1)^n/q^2 + \dots \; ,  $$
where $c_k = (2k)!/(k!(k+1)!)$ are the {\bf Catalan numbers}.
\end{lemma}

{\bf Remark.} This lemma 
was found empirically by determining 
the coefficients $c_k$ using high precision numerical methods using 1000 digit approximations of the algebraic numbers.
Looking up the resulting integer sequence revealed the Catalan connection. 
Even so I know now an elegant proof, it might be well known, but no reference so far has emerged. 
For $n=1$ for example, where $q=1, p/q=0$, the Catalan gold lemma gives the formula
$$  \sqrt{5} \alpha = \sum_{k=0}^{\infty} \frac{\left( \begin{array}{c} 2k \\ k \end{array} \right)} {(k+1) 5^k} \; . $$
The right hand side is $c(1/5) = 2/(1+\sqrt{1-4/5}) = \sqrt{5} (\sqrt{5}-1)/2$. 
The proof of the lemma has connections with Hilbert's 10th problem. The digging of Catalan Gold is only needed
to find accurate bounds for the interval in which $\sigma(y)$ converges and is left to the reader.

\begin{figure}
\scalebox{0.9}{\includegraphics{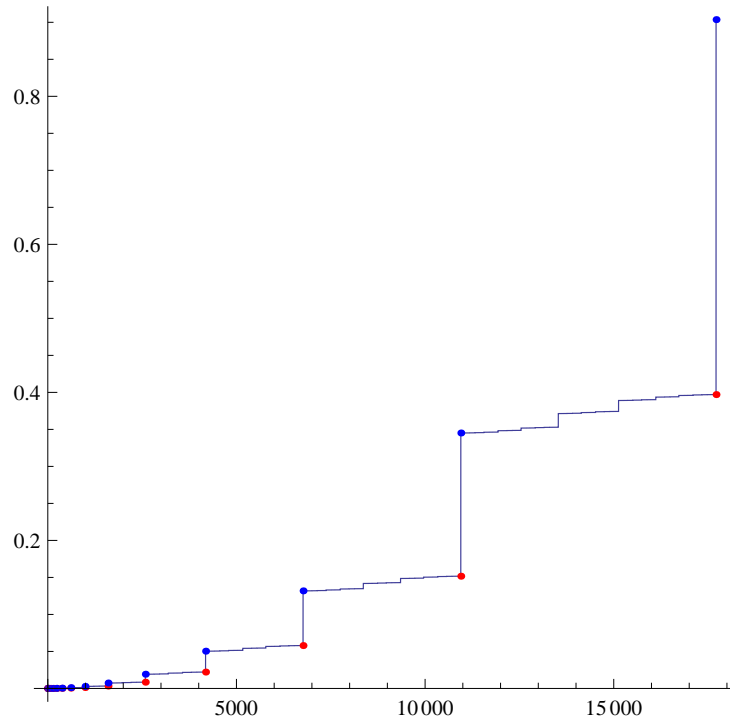}} \\
\caption{
The limiting function $S_{[x q_n]}/q_n^2$ of the $\csc^2$ Birkhoff sum
is monotone because the summands are positive. 
A modification of this graph has appeared in
\cite{KTp} as a discrepancy from the selfsimilarity of the difference graph. 
}
\end{figure}

\vspace{5mm}

\section{Remarks and outlook}

We briefly look at the Birkhoff sum for the antiderivative $G$ of $g$ which was 
studied in \cite{KTp}. The function has relations with Hecke's investigations \cite{hecke1922} because: 

\begin{propo}
$G(x)=\log(2-2 \cos(2\pi x))/2$ is the Hilbert transform of $H(x) = x-[x]-1/2$.
\end{propo}

\begin{proof}
Writing down some identities makes this clear:
$$ H(x) = \pi (x-[x]-1/2) = -\sum_{n=1}^{\infty} \frac{\sin(2\pi n x)}{n} = {\rm arg} (1-e^{2\pi i x})  \; , $$
$$ G(x) = \log(2-2 \cos(2\pi x))/2 = \log|2 \sin(\pi x)| 
        = -\sum_{n=1}^{\infty} \frac{\cos(2\pi nx)}{n} = {\rm log}|1-e^{2\pi i x}| \; .  $$
Note that while $G'(x)=\pi \cot(\pi x) = \pi g(x)$, the derivative of $H$ is only defined as a distribution.
There is no corresponding Hilbert dual result therefore for the $\cot$ function we were looking at. 
\end{proof}

{\bf Remark.}
The first identity for $H(x)$ is geometrically amusing: if a point $r(t)$ moves on a circle with uniform speed, 
then the direction vector from $r(t)$ to $r(0)$ changes with uniform speed. Since this can be done
from a second point also, it implies {\bf Thales result} about angles in a triangle where the base 
segment is fixed. 

As a direct consequence of functional analytic properties of the Hilbert transform, we get
a rather general Denjoy-Koksma type stability result. It confirms some numerical experiments in \cite{KL,KTp}:

\begin{coro}[Log boundedness of the Birkhoff sum] 
If $\alpha$ is Diophantine of constant type, then there exists a constant $C$ such that
for almost all $x$, we have 
$$   \frac{1}{\log(m)}  \sum_{k=1}^m G(x+k \alpha)  \leq C \; . $$
\end{coro}

\begin{proof}
The Hilbert transform is translational invariant. Therefore, the Birkhoff sum of the
conjugate functions remains the Hilbert transform of the original Birkhoff sum.
Since the Hilbert transform is a bounded operator on $L^2$, we have $||S_n(g,\alpha)||_2 \leq C ||S_n(h,\alpha)||_2$.
By Denjoy-Koksma applied to to function $h$ which has bounded variation, the later grows only logarithmic.
The former grows logarithmically too. 
\end{proof}

We still have to prove all the statements in \cite{KTp} but only outline here.
Lets denote by $S_n^G$ the Birkhoff sum for $G$, where the summation goes from $1$ to $q_n-1$.
We focus again on the situation, where $\alpha$ is the golden mean. The methods developed here
for $\cot$ will allow to show that the limit $S_n^G$ exists. 
The difference limit function is the sum of a self-similar
golden graph with different initial condition  and a correction which comes from
Abel summation and which is close to the graph of the $\csc$ Birkhoff limit.

\begin{theorem}
The Birkhoff difference functions
$f(x) = \lim_{n \to \infty} S_{[x q_n]}(p_{n+1}/q_{n+1})-S_{[x q_{n}]}(p_n/q_n)$
of $g(x) = \log(2-2 \cos(2\pi x))$ converge pointwise.
For the Riemann sum of $g(x) = \log(2-2\cos(x))$,
the limit $S_{q_n}$ converges and $\log(S_{q_n})/\log(q_n)$ stays bounded.
For the imaginary case,
$f_{2n}(x) \to c x^2, f_{2n-1}(x) \to -c x^2  \;  $ with $c=(1+\alpha^2) \pi/\sqrt{5}$.
Consequently, for the combination of real and imaginary part 
$g(x) = 2 \log(1-e^{ix})$ the difference limit exists.
The limiting function $f$ has the property that $f(\alpha) + \alpha^2 f$ is monotone.
\end{theorem}

\begin{proof}
The difference by a Rolle estimation, which is a Riemann sum for the
derivative $\cot$, but with a different starting point and a correction.
To estimate the difference function, we estimate the weighted sum
$$ f_n(x) = \sum_{j=1}^{[x q_n]} g(y+j \alpha) \; dx_j \; , $$
where $dx_j = p_k/q_k - p_{k+1}/q_{k+1}$.  \\
{\bf Abel summation} reduces this to the Birkhoff sum
$S_m= \sum_{k=1}^m g(k \alpha)$ with $g(x) = \cot(\pi x)$
for which we have shown the limit to exist. \\
The key is that for $g(x)=\cot(\pi x)$ and derivatives, the Rolle error function is bounded.
The monotonicity follows from the fact that $f(\alpha) + \alpha^2 f$ is a $\csc^2(x)$ Birkhoff sum, 
which is monotone.  The limit of $f_n$ is equal to the limit of
$\tilde{f}_n(x) = \sum_{j=1}^[x q_n] g(j \gamma_n) dx_j$
with $\gamma_n=(p_{n+1}/q_{n+1}+p_n/q_n)/2$.
Let $r_j$ denote the Rolle points. Then 
$$ f_n(x) = \sum_{j=1}^n g(r_j) \; dx_j $$
Abel summation allows to replace this by
$$ f_n(x) = S_{[x q_n]} \frac{1}{q_n} - \frac{1}{q_n^2} \sum_{k=1}^{[x q_n]} S_k  \; . $$
Because $S_k \leq k$, this converges. The first term is handled by the golden graph theorem.
The second part is just a Cesaro average.  \\
The imaginary part of $g$ is Hecke's Birkhoff sum $g_2(x) = 2\pi \{x\} = 2\pi ( x-[x]-1/2)$.
It leads to a limiting function satisfying $f(\alpha x) = \alpha^2 f(x)$ because $g''=0$.
The Abel sum gives
$$ f_n(x) = [x q_n] \frac{1}{q_n} - \frac{1}{q_n^2} \sum_{k=1}^{[x q_n]} 1 \; . $$
This explains why we have exact selfsimilarity and smoothness of the limiting difference.
Also this was experimentally first noted in \cite{KTp}.
\end{proof}

We could also modify the proof done for $\cot$. We have again a function $\tau(y)$ which is 
explicitly known. We write $\sigma(y)-\tau(y)$ as a sum and show convergence. 
We can do the summation even so we have no definite sign because we know the 
$t$ function and can use Abel summation to show that the limits exist.  \\
Why does $\sec(\pi x)=1/\cos(\pi x)$ give a different limiting golden graph than for $\cot$? 
The reason is a different pole at $x=\pi/2$ and is discontinuous at $0$. 
The trigonometric identity just mentioned shows however that we can 
reduce it as a sum of two $\cot$ Birkhoff sums with different starting points.  \\

\begin{figure}
\scalebox{0.9}{\includegraphics{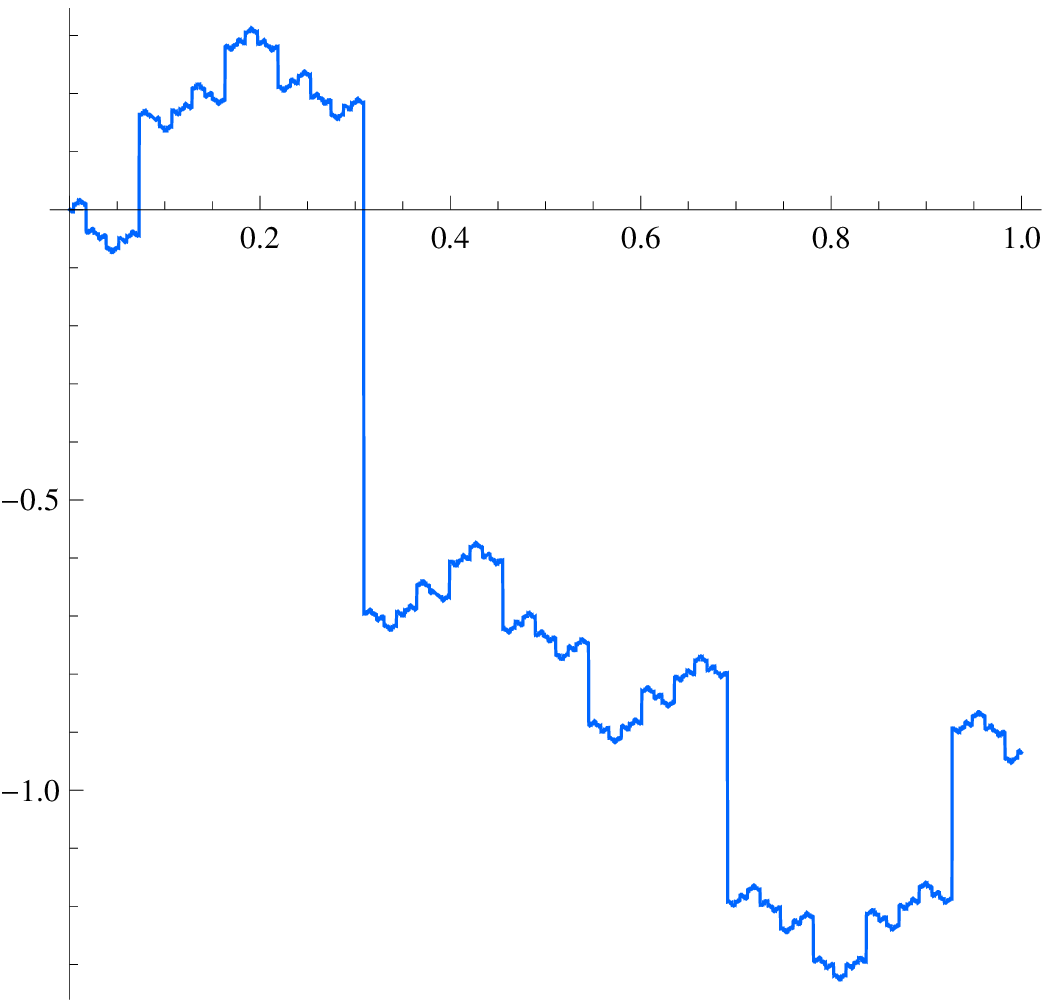}} \\
\caption{
The limiting Birkhoff sum function
$s'(x) = \lim_{n \to \infty} s_{6n}$ with the Birkhoff sum
$$  s'_n = \frac{1}{q_n} \sum_{k=1}^{[x q_n]} \sec(\pi k \alpha) $$
has period $6$. It is related to $R=\theta^2$ with
$R(\alpha) = 1+2 \sum_{k=1}^{\infty} \sec(\alpha)$, where ${\rm im}(\alpha)>0$. 
The attractor consists of 6 different limiting functions $\lim_{n \to \infty} s_{6n+k}$
with $k=0, \dots ,5$. We get the same attractor for $\cot(x+\pi/2)=-\tan(x)$ when
starting from $x=0$. 
}
\label{secsum}
\end{figure}

The function 
$R(\alpha) = \sum_{n \in Z} \sec(\pi n \alpha)$  satisfies $R(1/\alpha) = -\alpha R(\alpha)$
and $R(\alpha+2)=R(\alpha)$. It is a {\bf modular form} of weight $1$ and $R=\theta^2$, where
$\theta(\alpha) = \sum_{n \in Z} e^{i \pi n^2 \alpha} = \sum_{n \in Z} q^{n^2}$
which is a classical {\bf $\theta$ function}. \cite{C2} 
To understand the curlicue problem, finite sums of the later would be important. 


As an application we are so able to deal with the Birkhoff sum studied in \cite{KL,KTp}
which was historically first studied in \cite{HL46} who showed that 
$$ \frac{1}{n} \sum_{k=1}^n \log|\csc(k \alpha)| $$ 
converges to $\log(2)$ for almost all $\alpha$ implying that $\sum_{k=1}^n \log|2 \sin(k \alpha)| = o(n)$. 
We have seen in \cite{KL} that this is $O(\log(n)^2)$ for strongly Diophantine $\alpha$ and have seen
that for Diophantine $\alpha$ of constant type, the Birkhoff sum is $O(\log(n))$. \\

The difference limit function we have seen in \cite{KTp} is the sum of a selfsimilar
golden graph with different initial condition  and a correction which comes from
Abel summation and which is close to the graph of the $\csc$ Birkhoff limit of the $s'$.
The limiting behavior seen in \cite{KTp} can be explained with the $\cot$ series. \\

The following corollary is related to Lemmas 5 and 6 of Hardy-Littlewood's 1928 paper \cite{HL}, 
a paper which did not focus on the particular golden rotation numbers as we do. \\

\begin{coro}[Cosecant Birkhoff sum]
For small $|y|$, the Birkhoff sums
$S_{[x q_n]}(y)/q_n^2 = \sum_{k=1}^{[x q_n]} g'(y/q_n+k \alpha)$
of $g'(x) = - \pi \csc^2(\pi x)$ converge to a monotone function $s'_y(x)$
along even and odd subsequences. For $y=0$, the even and odd limits
agree and the graph satisfies $s'(\alpha x) = \alpha s'(x)$. 
\end{coro}

\begin{proof}
Differentiate the result for $\cot$. 
\end{proof}

What class of functions $g$ have the property that
$(1/q_n) \sum_{k=1}^{[x q_n]} g(k \alpha)$ converges along some subsequence?
We have seen that this is true for $f(x) = c \cot(\pi x) + h(x)$ where $h$
is continuous. This includes virtually all functions which have a single 
pole at $0$. A different type of pole appears with $\sec(\pi x)$ and also here, we have got 
limiting function, as seen in Figure~(\ref{secsum}). \\
For some other Diophantine numbers like quadratic irrationals, the story is similar. Instead
of a fixed point, we get a periodic attractor. In \cite{SU}, the expectation over $\alpha$ 
was taken, in which case one gets a limiting distribution and no subsequence has to be chosen. 
This suggests that for almost all $\alpha$, the renormalized random walks 
converge to an attractor with a limiting distribution. 
For Liouville $\alpha$, the situation is believed to be different and no limiting function 
$s(x)$ can exist. Indeed we can construct $\alpha$ close to rational numbers so that the Birkhoff
sum is unbounded.   \\
What about functions $g$ with more than one pole? An example is $g(x) + g(x+\beta)$ with $g(x) = \cot(\pi x)$.
If $\beta$ is rational, then a subsequence still converges and the graph is similar to a Birkhoff sum with one pole and 
quadratic irrational. For $\beta$ rationally independent from $1$ and $\alpha$ we do
not know what happens. The situation is even more difficult then than for non Diophantine
$\alpha$, where we still have the continued fraction expansion. The challenge is to determine,
how large the Birkhoff sum of $g(n \alpha) + g(n \alpha +\beta)$ can become. There is no reason why a Birkhoff 
limiting function should then exist. Indeed, experiments indicate that $s_m(x)$ explodes in general
with a speed which depends on arithmetic properties of $\alpha$ and $\beta$.

\bibliographystyle{plain}

\end{document}